\documentclass[12pt]{amsart}
\def\X{\mathcal X}

\def\Y{\mathcal Y}
\def\H{\mathcal H}
\def\B{\mathcal B}  
\def\iso{\cong}
\def\ptp{\hat{\otimes}}
\newcommand{\ig}{{\rm Im}\,}
 
\newcommand{\frechet}{Fr\'{e}chet}
\newcommand{\fr}{{\mathcal F} r}
\newcommand{\larray}{\left(\begin{array}{cc}}
\newcommand{\rarray}{\end{array}\right)}
\newtheorem{theorem}{\sc Theorem}[section]
\newtheorem{lemma}[theorem]{\sc Lemma}
\newtheorem{proposition}[theorem]{\sc Proposition}
\newtheorem{corollary}[theorem]{\sc Corollary}
\newtheorem{remark}[theorem]{\sc Remark}
\newtheorem{example}[theorem]{\sc Example}
\newtheorem{definition}[theorem]{\sc Definition}
\begin{document}

\title{Cyclic cohomology of certain nuclear Fr\'{e}chet and  $DF$ 
algebras}

\author{Zinaida A. Lykova}

\address{School of Mathematics and Statistics, University of Newcastle,\\
 Newcastle upon Tyne, NE1 7RU, UK~ {\rm (Z.A.Lykova@newcastle.ac.uk)}}

\date{22 August 2007}

\begin{abstract} We give explicit formulae for the continuous Hochschild and 
cyclic homology and cohomology of certain $\hat{\otimes}$-algebras. 
We use well-developed homological techniques together with some niceties of the theory of locally convex spaces to generalize the results known in the case of Banach algebras and their inverse limits to wider classes of 
topological algebras.
To this end we show that,  for a continuous morphism 
$ \varphi: \X\rightarrow \Y$ of complexes of complete nuclear $DF$-spaces, the isomorphism of cohomology groups
$H^n(\varphi):  H^n(\X) \rightarrow H^n(\Y)$
is automatically topological. The continuous cyclic-type
 homology and cohomology are described up to topological 
isomorphism for the following classes of biprojective 
$\hat{\otimes}$-algebras:  the tensor algebra $E \hat{\otimes} F$
 generated by the duality  $(E, F, \langle \cdot, \cdot \rangle)$
for nuclear \frechet\ spaces $E$ and $F$ or for nuclear 
$DF$-spaces $E$ and $F$; nuclear biprojective 
K\"{o}the algebras $\lambda(P)$ which are Fr\'echet spaces or
$DF$-spaces; the algebra of distributions $\mathcal{E}^*(G)$ 
on a compact Lie group $G$.\\

\end{abstract}

\thanks{I am indebted to the Isaac Newton Institute for Mathematical 
Sciences at Cambridge for hospitality and for generous financial 
support from the programme on Noncommutative Geometry while this work was 
carried out.}


\maketitle
	
\markboth{Z.~A.~Lykova}{Cyclic cohomology of nuclear Fr\'echet and  $DF$ 
algebras}

\section{Introduction}
Cyclic cohomology groups of topological algebras play an 
essential role in noncommutative geometry \cite{Co2}. 
There has been a number of papers addressing the calculation of
cyclic-type continuous homology and cohomology groups of some
Banach, $C^*$- and topological algebras; see, e.g.,
\cite{Co2,He1,Kh1,Ly3,Ly4,Mey,Wo1}. However, it remains 
 difficult to describe these groups explicitly for
many topological algebras. To compute the 
continuous Hochschild and cyclic cohomology groups  of 
\frechet\ algebras one has to deal with complexes
of complete $DF$-spaces. Here, in addition to presenting known homological techniques we also supply  technical enhancements that permit the
necessary generalization of results known in the case of Banach algebras and their inverse limits  to wider classes of topological algebras notably to those that occur in noncommutative geometry.

The category of Banach spaces has the useful property that
it is closed under passage to dual spaces. Fr\'echet spaces
do not have this property: the strong dual of a  Fr\'echet space is a
complete $DF$-space. $DF$-spaces have the awkward feature that 
their closed subspaces need not be  $DF$-spaces. However,  
closed subspaces of complete {\em  nuclear} $DF$-spaces are again 
 $DF$-spaces \cite[Proposition 5.1.7]{Pi}.

In Section 3 we use the strongest known results on the open
mapping theorem to give sufficient conditions on topological spaces
$E$ and $F$ to imply that any 
continuous linear operator $T$ from $E^*$ onto $F^*$ is open.
This allows us to prove the following results. 
In Lemma \ref{OpenMap-strong-dual-nuclear} we show that, 
 for a continuous morphism $ \varphi: \X
\rightarrow \Y$ of complexes of complete nuclear $DF$-spaces,
 the isomorphism of cohomology groups
$H^n(\varphi):  H^n(\X) \rightarrow H^n(\Y)$
is automatically topological.

 We use this fact to describe explicitly up to {\em topological} 
isomorphism the  continuous Hochschild and 
cyclic  cohomology groups of nuclear $\hat{\otimes}$-algebras $\mathcal{A}$
which are Fr\'echet spaces or  $DF$-spaces
and have trivial Hochschild homology ${\mathcal H}{\mathcal H}_n(\mathcal{A})$
for all $n \ge 1$ (Theorem \ref{A-simpl-trivial-Fr-DF}).
In Proposition \ref{A-simpl-trivial}, under the same condition on
${\mathcal H}{\mathcal H}_n(\mathcal{A})$, we give explicit formulae, 
 up to isomorphism of linear spaces, 
for continuous cyclic-type homology of $\mathcal{A}$ in a more general category of 
underlying spaces.

In Theorem \ref{cyclic-biproj-examples} the continuous cyclic-type
 homology and cohomology groups are described up to topological 
isomorphism for the following classes of biprojective 
$\hat{\otimes}$-algebras: the tensor algebra $E \hat{\otimes} F$
 generated by the duality  $(E, F, \langle \cdot, \cdot \rangle)$
for nuclear \frechet\ spaces  or for nuclear complete
$DF$-spaces $E$ and $F$; nuclear biprojective 
\frechet\ K\"{o}the algebras $\lambda(P)$;
 nuclear biprojective  K\"{o}the algebras $\lambda(P)^*$
which are $DF$-spaces;  the algebra of distributions $\mathcal{E}^*(G)$ 
 and  the algebra of smooth functions $\mathcal{E}(G)$ on a compact Lie 
group $G$. 
 
\section{Definitions and notation}
 We recall some notation and terminology used in homology 
and in the theory of topological algebras. Homological theory can be found in any 
relevant textbook, for instance, Loday \cite{Loday} for the pure algebraic
case and  Helemskii \cite{He0} for the continuous case.

Throughout the paper $\hat{\otimes}$ is the projective tensor product of 
complete locally convex spaces, by $ X^{\ptp n}$ we mean the $n$-fold projective tensor power $ X\ptp \dots \ptp X $ of $X$, and  ${\rm id}$
denotes the identity operator. 

We use the notation ${\mathcal Ban}$,  $\fr$  and $\mathcal{LCS}$
for the categories whose objects are Banach spaces, Fr\'{e}chet 
spaces and complete Hausdorff locally convex
spaces respectively, and whose morphisms in all cases are continuous
linear operators. 
For topological homology theory it is important to 
find a suitable category for the underlying spaces of the algebras 
and modules. In \cite{He0} Helemskii constructed homology theory for 
the following categories $\Phi$ of underlying spaces, for which he used the
notation $(\Phi, \hat{\otimes})$.

\begin{definition}\label{category}(\cite[Section II.5]{He0}) A {\em suitable category} for underlying spaces of the algebras and modules is
an arbitrary complete subcategory $\Phi$ of $\mathcal{LCS}$ having 
the following properties:

(i) if  $\Phi$ contains a space, it also contains all
those spaces topologically isomorphic to it;

(ii) if  $\Phi$ contains a space, it also contains 
any of its closed subspaces and the completion of any its Hausdorff
quotient spaces;

(iii)  $\Phi$ contains the direct sum and the
projective tensor product of any pair of its spaces;

(iv) $\Phi$ contains  ${\mathbf C}$. 
\end{definition}

Besides ${\mathcal Ban}$,  $\fr$  and $\mathcal{LCS}$ important
examples of suitable categories $\Phi$ are the categories of
complete nuclear spaces \cite[Proposition 50.1]{Tr}, 
 nuclear \frechet\ spaces and complete nuclear $DF$-spaces.
As to the above properties for the category of complete 
nuclear $DF$-spaces, recall the following results.
By \cite[Theorem 15.6.2]{Ja}, if $E$ and $F$ are 
complete $DF$-spaces, then $E \hat{\otimes} F$ is a 
 complete $DF$-space. By \cite[Proposition 5.1.7]{Pi},
a closed linear subspace of a complete nuclear $DF$-space is also
a complete nuclear $DF$-space. By \cite[Proposition 5.1.8]{Pi},
each quotient space of a complete nuclear $DF$-space by 
a closed linear subspace  is also a complete nuclear $DF$-space. 

By definition a $\hat{\otimes}$-algebra is a complete Hausdorff  
locally convex  algebra with jointly continuous multiplication. 
A left  $\hat{\otimes}$-module $X$ over a $\hat{\otimes}$-algebra ${\mathcal A}$ is a complete Hausdorff  locally convex 
space  $X$ together with the structure of a left ${\mathcal A}$-module such 
that the map 
${\mathcal A} \times X \to X$, $ (a,x) \mapsto a \cdot x $ is jointly continuous.
 For a  $\hat{\otimes}$-algebra 
${\mathcal A}$, $\hat{\otimes}_{\mathcal A}$ is the projective tensor product
of   left and right  ${\mathcal A}$-$\hat{\otimes}$-modules (see \cite{GrA},
\cite[II.4.1]{He0}).
The category of  left ${\mathcal A}$-$\hat{\otimes}$-modules is denoted by 
${\mathcal A}$-{\rm mod}  and the category of  ${\mathcal A}$-$\hat{\otimes}$-bimodules 
is denoted by  ${\mathcal A}$-{\rm mod}-${\mathcal A}$. 

 Let ${\mathcal K}$ be one of the above categories. 
 A {\it  chain complex} ${\mathcal X}_{\sim} $ in the category  ${\mathcal K}$
 is a  sequence of $ X_n \in {\mathcal K}$ and  morphisms $d_n $ 
$$ \dots \leftarrow X_n \stackrel {d_n} {\leftarrow}  X_{n+1} 
\stackrel {d_{n+1}} {\leftarrow} X_{n+2}
 \leftarrow  \dots $$
such that $d_n \circ d_{n+1} = 0 $ for every $n$.
The {\it homology groups of ${\mathcal X}_{\sim} $} are defined by 
$$ H_n({\mathcal X}_{\sim}) = {\rm Ker}~ d_{n-1}/{\rm Im}~ d_n.$$
  A continuous morphism of chain complexes 
${\psi}_{\sim} : {\mathcal X}_{\sim} \rightarrow {\mathcal P}_{\sim}$
induces a continuous linear operator 
$H_n({\psi}_{\sim}) : H_n({\mathcal X}_{\sim} )
 \rightarrow H_n({\mathcal P}_{\sim})$ \cite[Definition 0.4.22]{He3}.

If $E$ is  a topological vector space $E^*$ denotes its dual space of
continuous linear functionals.
 Throughout the paper, $E^*$ will always be equipped with the strong topology
unless otherwise stated. The  {\it strong topology} is defined  on $E^*$ by taking as a basis of 
neighbourhoods of $0$ the family of polars $V^0$ of all bounded subsets $V$ of $E$; see \cite[II.19.2]{Tr}.

 For any $\hat{\otimes}$-algebra ${\mathcal A}$, not necessarily unital, 
${\mathcal A}_+$ is the  $\hat{\otimes}$-algebra obtained by 
adjoining an identity to ${\mathcal A}$.
 For a $\hat{\otimes}$-algebra ${\mathcal A}$, the algebra 
${\mathcal A}^e = {\mathcal A}_+ \hat{\otimes} {\mathcal A}_+^{op}$
is called the {\it 
 enveloping algebra of ${\mathcal A}$}, where ${\mathcal A}_+^{op}$ is the {\it  opposite algebra}
of ${\mathcal A}_+$ with  multiplication $a \cdot b = ba$.

 A complex of  ${\mathcal A}$-$\hat{\otimes}$-modules and their morphisms
is called {\it  admissible} if it 
splits as a complex in  $\mathcal{LCS}$ \cite[III.1.11]{He0}.
 A module $Y \in {\mathcal A}$-{\rm mod} is called 
{\it flat } if for any admissible complex ${\X}$ of right 
 ${\mathcal A}$-$\hat{\otimes}$-modules the complex 
${\X} \hat{\otimes}_{{\mathcal A}} Y$ is exact. 
A module $Y \in  {\mathcal A}$-{\rm mod}-${\mathcal A}$ is called 
{\it flat } if for any admissible complex ${\X}$ of  
 ${\mathcal A}$-$\hat{\otimes}$-bimodules the complex 
${\X} \hat{\otimes}_{{\mathcal A}^e} Y$ is exact.
 For $Y, X  \in {\mathcal A}$-{\rm mod}-${\mathcal A}$, we shall denote 
by  $ {\rm Tor}_{n}^{{\mathcal A}^e}(X, Y)$ the $n$th homology 
of the complex  $ X \ptp_{{\mathcal A}^e}{\mathcal P}$, where 
$0 \leftarrow Y \leftarrow {\mathcal P}$ is a  projective resolution 
of $Y$ in ${\mathcal A}$-{\rm mod}-${\mathcal A}$, \cite[Definition III.4.23]{He0}.

It is well known that the strong dual of a Fr\'{e}chet space is a 
complete $DF$-space and that nuclear Fr\'{e}chet spaces and complete
nuclear $DF$-spaces are reflexive
\cite[Theorem 4.4.12]{Pi}. Moreover, the correspondence 
$E \leftrightarrow E^*$  establishes a one-to-one relation between nuclear 
Fr\'{e}chet spaces and complete nuclear $DF$-spaces 
\cite[Theorem 4.4.13]{Pi}. $DF$-spaces were introduced 
by A. Grothendieck in \cite{GrA1}.

Further we shall need the following technical result which extends a result of Johnson for the Banach case \cite[Corollary 1.3]{BEJ1}.

\begin{proposition}\label{n>N-homology-dual-cohomology}
Let $(\mathcal X, d)$ be a chain complex of \\
{\rm (a)} Fr\'{e}chet spaces and continuous linear operators, or\\
{\rm (b)} complete nuclear $DF$-spaces and continuous linear 
operators,\\
and let $N \in {\mathbf N}$. Then the following statements are equivalent:

\noindent{\rm (i)} $H_n(\mathcal X, d)= \{0\}\;$ for all $n \ge N$
and $H_{N-1}(\mathcal X, d)$ is Hausdorff; 

\noindent{\rm (ii)} $H^n({\mathcal X}^*, d^*)= \{0\}\;$ for all $n \ge N.$

\end{proposition}

\begin{proof} ~ Recall that $H_n(\mathcal X, d)= 
{\rm Ker}~d_{n-1}/{\rm Im}~d_n$
and $H^n({\mathcal X}^*, d^*)= {\rm Ker}~d^*_{n}/{\rm Im}~d^*_{n-1}.$
Let $L$ be the closure of ${\rm Im}~ d_{N-1}$  in $X_{N-1}$. 
Consider the
following commutative diagram
\begin{equation} 
\begin{array}{ccccccccccc}
\label{n>N-homology} 
0 & \leftarrow & L &\stackrel{j} {\longleftarrow}& X_{N}
 & \stackrel{d_{N}} {\longleftarrow}  & X_{N+1}
& \stackrel{d_{N+1}} {\longleftarrow}  & \dots\\
~ & ~ & \downarrow\vcenter{\rlap{$\scriptstyle{i}~~$}} &  
\swarrow\vcenter{\rlap{$\scriptstyle{d_{N-1}}~~$}} ~ &  ~ &~ &~ &~  \\
~ & ~ & X_{N-1}& ~ & ~ &~ & ~ &~ &\\ 
\end{array} 
\end{equation}
in which $i$ is the natural inclusion and $j$ is a corestriction of 
$d_{N-1}$. The dual commutative diagram  is the
following
\begin{equation} 
\begin{array}{ccccccccccc} 
\label{n>N-cohomology}
0 & \rightarrow & L^* &\stackrel{j^*} {\longrightarrow}& X_{N}^*
 & \stackrel{d_{N}^*} {\longrightarrow}  & X_{N+1}^*
& \stackrel{d_{N+1}^*} {\longrightarrow}  & \dots\\
~ & ~ & \uparrow\vcenter{\rlap{$\scriptstyle{i^*}~~$}} &  
\nearrow\vcenter{\rlap{$\scriptstyle{d_{N-1}^*}~~$}} ~ &  ~ &~ &~ &~  \\
~ & ~ & X_{N-1}^*& ~ & ~ &~ & ~ &~ &\\ 
\end{array} 
\end{equation}
It is clear that  $H_{N-1}(\mathcal X, d)$ is Hausdorff if and only if
$j$ is surjective. Since $i$ is injective, condition (i) is equivalent
to the exactness of  diagram (\ref{n>N-homology}). On the other hand,
by the Hahn-Banach theorem,  $i^*$ is surjective. Thus condition (ii) 
is equivalent to the exactness of  diagram (\ref{n>N-cohomology}).

In the case of  Fr\'{e}chet spaces, by \cite[Lemma 2.3]{Ly4}, the 
exactness of the complex  (\ref{n>N-homology}) is  equivalent to 
the exactness of the complex  (\ref{n>N-cohomology}).

In the case of  complete nuclear 
$DF$-spaces, by \cite[Proposition 5.1.7]{Pi}, $L$ is the strong dual of 
a nuclear Fr\'{e}chet space. By \cite[Theorem 4.4.12]{Pi}, complete
nuclear $DF$-spaces are reflexive, and therefore the complex  
(\ref{n>N-homology}) is the dual of the complex  
(\ref{n>N-cohomology}) of nuclear Fr\'{e}chet spaces and continuous 
linear operators. By \cite[Lemma 2.3]{Ly4}, the 
exactness of the complex  (\ref{n>N-homology}) is  equivalent to 
the exactness of the complex  (\ref{n>N-cohomology}).
The proposition is proved.
\end{proof}

\section{The open mapping theorem in complete nuclear $DF$-spaces}

It is known that there exist closed linear subspaces of $DF$-spaces 
that are not $DF$-spaces. For nuclear spaces, however, we have the
following.

\begin{lemma}\label{closed_suspace_DFspaces} \cite[Proposition 5.1.7]{Pi}
Each closed linear subspace $F$ of the strong dual $E^*$ of
a nuclear Fr\'echet space $E$ is also the strong dual of a nuclear 
Fr\'echet space.
\end{lemma}

In a locally convex space a subset is called a {\it barrel} if it 
is absolutely convex, absorbent and closed. Every locally
convex space has a neighbourhood base consisting of barrels. 
A locally convex space is called  a {\it barrelled} space if every 
barrel is a neighbourhood \cite{RR}. By \cite[Theorem IV.1.2]{RR}, every 
Fr\'echet space is barrelled. By \cite[Corollary IV.3.1]{RR}, a Hausdorff locally convex space is reflexive if and only if it is barrelled and every bounded set is contained in a weakly compact set. Thus the strong dual of a 
nuclear Fr\'echet space is barrelled.
For a generalization of the open mapping theorem to locally convex spaces,
V. Pt\'{a}k introduced the notion of $B$-completeness in \cite{Pt}.
A subspace $Q$ of  $E^*$ is said to be {\it almost closed} if, for each neighbourhood $U$ of $0$ in $E$, $Q \cap U^0$ is closed in the relative 
weak* topology $\sigma(E^*, E)$ on $U^0$.
A locally convex space $E$ is said to be {\it $B$-complete} or {\it fully complete} if each almost closed subspace of $E^*$ is closed in the weak* topology $\sigma (E^*, E)$. 

\begin{theorem}\label{open-map-Ptak} \cite{Pt}. 
Let $E$ be a $B$-complete locally convex
space and $F$ be a  barrelled locally convex space. Then a continuous 
 linear operator  $f$ of $E$ onto $F$ is open.
\end{theorem}

Recall \cite[Theorem 4.1.1]{Hu} that a locally convex space $E$ 
is $B$-complete if and only if each linear continuous  and almost open 
mapping $f$ of $E$ onto any locally convex space $F$ is open. 
By  \cite[Proposition 4.1.3]{Hu},  every Fr\'echet space  is $B$-complete.

\begin{theorem}\label{open-map-strong-dual}  Let $E$ be a semi-reflexive 
metrizable barrelled space, 
 $F$ be a  Hausdorff reflexive locally convex space  and let
$E^*$ and $F^*$ be the strong duals of  $E$ and $F$ respectively.
Then a continuous linear operator $T$ of $E^*$ onto $F^*$ is open.
\end{theorem}
\begin{proof} By \cite[Theorem 6.5.10]{Hu} and by 
\cite[Corollary 6.2.1]{Hu}, the strong dual $E^*$ of a semi-reflexive
metrizable barrelled space $E$ is $B$-complete.
By \cite[Corollary IV.3.2]{RR}, if a Hausdorff locally convex space is 
reflexive, so is its dual under the strong topology.
By \cite[Corollary IV.3.1]{RR}, a Hausdorff reflexive locally convex space 
is barrelled. Hence $F^*$ is  a  barrelled locally convex space.
Therefore, by Theorem \ref{open-map-Ptak},  $T$ is open.
\end{proof}

\begin{corollary}\label{open-map-strong-dual-FM} Let $E$ and $F$
be  nuclear Fr\'echet spaces and let
$E^*$ and $F^*$ be the strong duals of  $E$ and $F$ respectively.
Then a  continuous linear operator  $T$ of $E^*$ onto $F^*$ is open.
\end{corollary}

For a continuous morphism of chain complexes 
${\psi}_{\sim} : {\mathcal X}_{\sim} \rightarrow {\mathcal P}_{\sim}$
in $\fr $, a surjective map 
$ H_n(\varphi):  H_n(\X) \rightarrow H_n(\Y)$ is automatically
open, see \cite[Lemma 0.5.9]{He0}.
To get the corresponding result for dual complexes of Fr\'echet spaces
one has to assume nuclearity.

\begin{lemma}\label{OpenMap-strong-dual} Let  $(\X, d_{\X}) $ and $(\Y,
d_{\Y})$ be  chain complexes  of nuclear Fr\'echet spaces and continuous
linear operators and let $ (\X^*,d_{\X}^* )$  and $(\Y^*,
d_{\Y}^*)$ be their strong dual complexes.
Let $ \varphi: \X^*
\rightarrow \Y^*$ be a continuous morphism of complexes.
Suppose that
$$\varphi_* = H^n(\varphi): 
 H^n(\X^*,d_{\X}^* ) \rightarrow H^n (\Y^*,d_{\Y}^*)$$
is surjective. Then $\varphi_*$ is open.
\end{lemma}
\begin{proof} 
Let 
$~\sigma_{\Y^*}: {\rm Ker}~ (d_{\Y}^*)_{n} \to H^n(\Y^*,d_{\Y}^*)~$ be the
quotient map. Consider the map 
$$\psi: {\rm Ker}~ (d_{\X}^*)_{n} \oplus Y_{n-1}^* \rightarrow {\rm Ker}~ (d_{\Y}^*)_{n}
\subset Y_{n}^*$$
given by $ (x,y) \mapsto \varphi_n(x) + (d_{\Y}^*)_{n-1}(y).$

By Lemma \ref{closed_suspace_DFspaces}, ${\rm Ker}~  (d_{\X}^*)_{n}$  and  
${\rm Ker}~ (d_{\Y}^*)_{n}$ are 
the strong duals of nuclear Fr\'echet spaces and hence are barrelled. 
By \cite[Theorem 6.5.10]{Hu} and  
\cite[Corollary 6.2.1]{Hu}, the strong dual of a semi-reflexive
metrizable barrelled space is $B$-complete. Thus 
${\rm Ker}~ (d_{\X}^*)_{n}$, $ Y_{n-1}^* $  and 
$$
{\rm Ker}~ (d_{\X}^*)_{n} \oplus Y_{n-1}^* \iso 
[({\rm Ker}~ (d_{\X}^*)_{n})^* \oplus Y_{n-1}]^* 
$$
are  $B$-complete.
By assumption $\varphi_*$ maps $H^n(\X^*,d_{\X}^* )$ onto 
$H^n(\Y^*,d_{\Y}^*)$, which implies that $\psi$
is a surjective linear continuous operator from the $B$-complete 
locally convex space
${\rm Ker}~ (d_{\X}^*)_{n} \oplus Y_{n-1}^* $ to the barrelled
locally convex space $ {\rm Ker}~ (d_{\Y}^*)_{n}$.
Therefore, by Theorem \ref{open-map-Ptak},
$\psi$ is open. Consider the diagram 
\begin{equation} 
\begin{array}{ccccccccc} 
{\rm Ker}~ (d_{\X}^*)_{n} \oplus Y_{n-1}^*  & \stackrel{j} {\rightarrow}  & 
{\rm Ker}~ (d_{\X}^*)_{n} &
\stackrel{\sigma_{\X^*}} {\rightarrow} &  H^n(\X^*,d_{\X}^* )\\
\downarrow\vcenter{\rlap{$\scriptstyle{\psi}~~$}} & ~ & ~ & ~ &
\downarrow\vcenter{\rlap{$\scriptstyle{\varphi_*}~~$}}\\ {\rm Ker}~
(d_{\Y}^*)_{n} & ~ &
\stackrel{\sigma_{\Y^*}} {\longrightarrow} &~ & H^n(\Y^*,d_{\Y}^*)\\ 
\end{array} 
\end{equation}
in which $j$ is a projection onto a direct summand and $\sigma_{\X^*}$ 
and $\sigma_{\Y^*}$ are the natural quotient maps.
Obviously this diagram is commutative. Note that the projection
$j$ and quotient maps $\sigma_{\X^*}$,  $\sigma_{\Y^*}$ are  open. 
As $\psi$ is also an open map, so is 
$\sigma_{\Y^*} \circ \psi =\varphi_* \circ
\sigma_{\X^*} \circ j$. Since $\sigma_{\X^*} \circ j$ is continuous, 
 $\varphi_*$ is open.
\end{proof}

\begin{corollary}\label{OpenMap-strong-dual-nuclear} Let 
 $(\X, d_{\X}) $ and $(\Y, d_{\Y})$ be  cochain complexes of 
complete nuclear $DF$-spaces and continuous
linear operators, and let $ \varphi: \X \rightarrow \Y$ be a 
continuous morphism of complexes.
Suppose that $ \varphi_* = H^n(\varphi): 
 H^n(\X,d_{\X}) \rightarrow H^n (\Y,d_{\Y})$
is surjective. Then $\varphi_*$ is open.
\end{corollary}
\begin{proof} By \cite[Theorem 4.4.13]{Pi}, 
$(\X, d_{\X}) $ and $(\Y, d_{\Y})$ are strong duals of 
chain complexes  $ (\X^*,d_{\X}^* )$  and $(\Y^*,d_{\Y}^*)$
of  nuclear Fr\'{e}chet spaces  and continuous operators.
The result follows from Lemma \ref{OpenMap-strong-dual}.
\end{proof}

\section{Cyclic and Hochschild cohomology of some  
$\hat{\otimes}$-algebras}
One can consult the books by Loday \cite{Loday} or Connes 
\cite{Co2} on cyclic-type homological theory.

Let ${\mathcal A}$ be a $\hat{\otimes}$-algebra and let $X$ be an 
${\mathcal A}$-$\hat{\otimes}$-bimodule.  We assume here that the category 
of underlying spaces $\Phi$ has the properties from Definition \ref{category}.
  Let us recall the definition of the standard homological
chain complex ${\mathcal C}_{\sim} ({\mathcal A}, X)$. For $n\ge 0$, let $C_n({\mathcal A}, X)$
denote the projective tensor product 
 $X \ptp {\mathcal A}^{{\ptp}^n}$. 
The elements of $C_n({\mathcal A}, X)$ are called
{\em $n$-chains}. Let the differential $d_n : C_{n+1} \to C_n$ be given by
$$d_n(x \otimes  a_1 \otimes \ldots \otimes  a_ {n+1})
= x \cdot a_1 \otimes \ldots \otimes  a_ {n+1}+$$
$$  \sum_{k=1}^{n} (-1)^k (x \otimes a_1 \otimes \ldots \otimes a_k
a_{k+1} \otimes \ldots \otimes a_ {n+1})
 +(-1)^{n+1}(a_ {n+1} \cdot x  \otimes a_1 \otimes \ldots \otimes 
a_{n})$$
with $d_{-1}$ the null map. The homology groups of  this complex
$H_n({\mathcal C}_{\sim} ({\mathcal A}, X))$ are called  the
{\it continuous Hochschild homology groups of  ${\mathcal A}$ 
with coefficients in $X$} and
denoted by $\H_n({\mathcal A}, X)$ \cite[Definition II.5.28]{He0}. 
We also consider the cohomology groups 
$H^n(({\mathcal C}_{\sim} ({\mathcal A}, X))^*)$
of the dual complex $({\mathcal C}_{\sim} ({\mathcal A}, X))^*$ 
with the strong dual topology. For Banach algebras ${\mathcal A}$,
$H^n(({\mathcal C}_{\sim} ({\mathcal A}, X))^*)$ is topologically isomorphic 
to the Hochschild cohomology  $\H^n({\mathcal A}, X^*)$ of ${\mathcal A}$
with coefficients in the dual ${\mathcal A}$-bimodule $X^*$ 
\cite[Definition I.3.2 and Proposition II.5.27]{He0}. 
The {\it  weak bidimension} of a Fr\'{e}chet algebra ${\mathcal A}$ is
$${\rm db}_w {\mathcal A} = \inf \{ n: 
H^{n+1}({\mathcal C}_{\sim} ({\mathcal A}, X)^*) = \{ 0 \} \; 
{\rm for \; all \; Fr\acute{e}chet} \;{\mathcal A}{\rm -bimodules} \;  X \}.
$$

The  {\em continuous bar} and {\em \lq naive' Hochschild
homology of} a  $\hat{\otimes}$-algebra ${\mathcal A}$ are 
defined respectively as 
$$
{\mathcal H}^{bar}_*({\mathcal A}) = H_*({\mathcal C}({\mathcal A}), b') 
\;\; {\rm and}\;\; 
{\mathcal H}^{naive}_*({\mathcal A})= H_*({\mathcal C} ({\mathcal A}), b),$$
where  ${\mathcal C}_n({\mathcal A}) =  {\mathcal A}^{\hat{\otimes}(n+1)}$,  and the differentials
$b$, $b'$ are given by 
$$
b'(a_0 \otimes \dots \otimes a_n) = \sum_{i=0}^{n-1}(-1)^i (a_0 \otimes \dots
\otimes a_i a_{i+1}\otimes \dots \otimes a_n) \; {\rm and} \;
$$
$$
b(a_0 \otimes \dots \otimes a_n) = 
b'(a_0 \otimes \dots \otimes a_n) + 
(-1)^n(a_n a_0 \otimes \dots \otimes a_{n-1}).
$$
Note that ${\mathcal H}^{naive}_*({\mathcal A})$ is just another way of writing 
${\mathcal H}_*({\mathcal A}, {\mathcal A})$,
the continuous  homology of ${\mathcal A}$ with coefficients in ${\mathcal A},$ as described in 
 \cite{He0,BEJ1}.

There is a powerful method based on mixed complexes for the study of 
the  cyclic-type homology groups; see papers by C. Kassel \cite{Kas1},
 J. Cuntz and D. Quillen \cite{CQ2} and J. Cuntz \cite{Cu}. We shall 
present this method for
the category $\mathcal{LCS}$ of locally convex spaces and 
continuous linear operators; see \cite{BL} for the category
of Fr\'{e}chet spaces. A {\em mixed complex} 
$(\mathcal{M}, b, B)$ in
the category $\mathcal{LCS}$ is a family  $\mathcal{M} = \{ M_n\}_{n\ge 0}$ 
of locally convex spaces $M_n$ equipped with  
{\em continuous} linear operators $b_n: M_n \rightarrow M_{n-1}$ and  $B_n: M_n \rightarrow
M_{n+1}$, which  satisfy the identities 
$b^2 = bB + Bb = B^2 = 0$. We assume that in
degree zero the  differential $b$ is identically equal to zero.
We arrange the mixed complex $(\mathcal{M}, b, B)$ in the double complex
\begin{equation}\label{familiar-double-complex}
\begin{array}{ccccccc}
\dots & ~ & \dots & ~ & \dots & ~ & \dots\\
b \downarrow & & b \downarrow & &b \downarrow \\
M_2 & \stackrel {B}{\leftarrow} & M_1 & \stackrel {B}{\leftarrow}
 &  M_0 \\
b \downarrow && b \downarrow\\
M_1 & \stackrel {B}{\leftarrow} & M_0 \\
b \downarrow\\
M_0\\
\end{array}
\end{equation}
There are three types of homology theory that can be naturally associated 
with a mixed complex. The {\it  Hochschild 
homology}
 $H^b_*(\mathcal{M})$ of $(\mathcal{M}, b , B)$  is the homology of 
the chain complex 
$(\mathcal{M}, b)$, that is,  
$$
H^b_n(\mathcal{M}) = H_n(\mathcal{M}, b)
 = {\rm Ker}\;\{ b_n: M_n \rightarrow M_{n-1}\}/ \ig \{ b_{n+1}: 
M_{n+1} \rightarrow M_n\}.
$$
To define the cyclic homology of $(\mathcal{M} , b, B)$, 
let us denote by $\mathcal{B}_c \mathcal{M}$ the total complex of the above
double complex, that is,
$$
\dots \rightarrow (\mathcal{B}_c\mathcal{M})_{n}
\stackrel {b+B} {\rightarrow}  (\mathcal{B}_c \mathcal{M})_{n-1}
\rightarrow \dots
\stackrel {b+B} {\rightarrow} 
   (\mathcal{B}_c\mathcal{M})_{0} \rightarrow 0, 
   $$ 
where the spaces
$$
(\mathcal{B}_c \mathcal{M})_0 = M_0, \;  \dots,\;
(\mathcal{B}_c \mathcal{M})_{2k -1} = M_1 \oplus M_3 \oplus \dots \oplus M_{2k -1}
$$
and
$$
 (\mathcal{B}_c \mathcal{M})_{2k} = M_0 \oplus M_2 \oplus \dots \oplus 
M_{2k}
$$
are equipped with the product
topology, and the continuous linear operators $b+B$ are defined 
by 
$$
(b+B)(y_0, \dots , y_{2k}) = 
( by_2 + B y_0, \dots, by_{2k}+ By_{2k-2} )\;
$$
and 
$$\; (b+B)(y_1, \dots , y_{2k+1}) = 
( by_1, \dots, by_{2k+1}+ By_{2k-1} ).
$$
The {\em cyclic homology} of $(\mathcal{M}, b, B)$
is  defined to be
$
H_*(\mathcal{B}_c \mathcal{M}, b+B).
$
It is denoted by $H^c_*(\mathcal{M}, b, B ).$

The periodic cyclic homology of $(\mathcal{M}, b, B)$ is defined  
in terms of the complex
$$
\dots \rightarrow (\mathcal{B}_p\mathcal{M})_{ev} \stackrel {b+B} {\rightarrow}
 (\mathcal{B}_p\mathcal{M})_{odd} \stackrel {b+B} {\rightarrow}
(\mathcal{B}_p\mathcal{M})_{ev} \stackrel {b+B} {\rightarrow} 
(\mathcal{B}_p\mathcal{M})_{odd}
\rightarrow \dots,
$$
where even/odd chains are elements
of the product spaces
$$
(\mathcal{B}_p \mathcal{M})_{ev} = \prod_{n\geq 0} M_{2n}\;\;  
{\rm and} \;\;  (\mathcal{B}_p \mathcal{M})_{odd} 
 = \prod_{n\geq 0} M_{2n +1}, 
$$
respectively. The spaces $(\mathcal{B}_p \mathcal{M})_{ev/odd}$  are 
locally convex
spaces with respect to the product topology \cite [Section 18.3.(5)]{Ko1}. 
The continuous differential $b+B$ is
defined as an obvious extension of the above.
The {\it periodic cyclic homology} of  
 $( \mathcal{M} , b, B)$ is   $
H^p_\nu (\mathcal{M}, b, B ) = H_\nu (\mathcal{B}_p \mathcal{M} , b+ B)
$, where $\nu \in {\bf Z}/2{\bf Z}$.

There are also three types of cyclic {\it cohomology} theory associated 
with the mixed complex,
obtained when one replaces the chain complex of locally convex
 spaces by its dual complex of strong dual spaces.
For example, the {\it cyclic  cohomology} associated 
with the mixed complex $(\mathcal{M}, b, B )$ is defined to be 
the cohomology of the dual complex $((\mathcal{B}_c \mathcal{M})^*, b^* +B^*)$
of strong dual spaces and dual operators; it is denoted by 
$H_c^* (\mathcal{M}^*, b^*, B^*).$

Consider the mixed complex 
$(\bar{\Omega} {\mathcal A}_+, \tilde{b}, \tilde{B})$, where
$\bar{\Omega}^n {\mathcal A}_+ = {\mathcal A}^{\hat{\otimes}(n+1)} \oplus 
{\mathcal A}^{\hat{\otimes} n}$ and
$$
\tilde{b} = \larray b & 1-\lambda \\ 0 & -b'\rarray; \;\;\;
\tilde{B} = \larray 0 & 0 \\ N & 0 \rarray 
$$
where $\lambda(a_1 \otimes \dots \otimes a_{n}) = 
(-1)^{n-1}(a_n \otimes a_1 \otimes \dots \otimes a_{n-1})$ and
$ N = {\rm id} + \lambda + \dots + \lambda^{n-1}$  \cite[1.4.5]{Loday}.  
The {\it continuous Hochschild homology of ${\mathcal A}$}, the {\it continuous cyclic
homology of ${\mathcal A}$ }  and the {\it  continuous 
periodic cyclic homology of ${\mathcal A}$ } are
defined by 
$$
{\mathcal H}{\mathcal H}_*({\mathcal A}) = H^b_*(\bar{\Omega} {\mathcal A}_+, \tilde{b}, \tilde{B}),
\; \;
{\mathcal H}{\mathcal C}_*({\mathcal A}) = H^c_*(\bar{\Omega} {\mathcal A}_+, \tilde{b}, \tilde{B})
\; \;{\rm and} $$
$$
{\mathcal H}{\mathcal P}_*({\mathcal A}) = H^p_*(\bar{\Omega} {\mathcal A}_+, \tilde{b}, \tilde{B})
$$
where $H^b_*$, $H^c_*$ and $H^p_*$ are 
Hochschild homology, cyclic homology and periodic cyclic homology
of the mixed complex $(\bar{\Omega} {\mathcal A}_+, \tilde{b}, \tilde{B})$ in
$\mathcal{LCS}$, see \cite{Ly3}.

There is also a {\it cyclic  cohomology} theory associated 
with a complete locally convex algebra ${\mathcal A},$
obtained when one replaces the chain complexes of ${\mathcal A}$ by their dual 
complexes of strong dual spaces.

\begin{lemma}\label{topology_of_complexes}  {\rm (i)} Let ${\mathcal A}$ be a [nuclear] Fr\'{e}chet algebra.
Then the following complexes 
$({\mathcal C} ({\mathcal A}), b)$, $(\bar{\Omega} {\mathcal A}_+, \tilde{b})$,
$(\mathcal{B}_c \bar{\Omega} {\mathcal A}_+, \tilde{b} +\tilde{B})$  and
$(\mathcal{B}_p \bar{\Omega} {\mathcal A}_+, \tilde{b} +\tilde{B})$ 
are complexes of [nuclear] Fr\'{e}chet spaces and continuous linear operators.

{\rm (ii)} Let ${\mathcal A}$ be a [nuclear] $\hat{\otimes}$-algebra which is a  $DF$-space. Then the following complexes 
$({\mathcal C} ({\mathcal A}), b)$, $(\bar{\Omega} {\mathcal A}_+, \tilde{b})$, and $(\mathcal{B}_c \bar{\Omega} {\mathcal A}_+, \tilde{b} +\tilde{B})$ 
are complexes of [nuclear] complete $DF$-spaces and continuous linear operators, and 
$(\mathcal{B}_p \bar{\Omega} {\mathcal A}_+, \tilde{b} +\tilde{B})$ is a
 complex of [nuclear] complete locally convex spaces and continuous linear operators, but it is not a $DF$-space in general.
\end{lemma}

\begin{proof} 
It is well known that  Fr\'{e}chet spaces are closed under countable cartesian products and projective tensor product \cite{Tr}; nuclear locally convex spaces are closed under cartesian products, countable direct sums and projective tensor product \cite[Corollary 21.2.3]{Ja}; complete $DF$-spaces are closed under countable direct sums, projective tensor product, but not under infinite cartesian products \cite[Theorem 12.4.8 and Theorem 15.6.2]{Ja}.
\end{proof}


Propositions \ref{n>N-HHn=0-homology} and \ref{n>N-HHn=0-cohomology} below are proved by the author in \cite{Ly3,Ly4} and show the equivalence between the  continuous cyclic (co)homology of $A$   and the  continuous periodic cyclic (co)homology of $A$  when $A$ has trivial
continuous Hochschild (co)homology $HH_n(A)$ for all $n \ge N$ for some
integer $N$. Here we add in these statements certain  topological conditions on the algebra which allow us to show that isomorphisms of 
(co)homology groups are automatically topological.

\begin{proposition}\label{n>N-HHn=0-homology}{\em \cite[Proposition 3.2]{Ly3}}
Let $A$ be a  complete locally convex  algebra.
Then, for  any even integer $N$, say
$N=2K$,  and the following assertions, we have 
$ {\rm (i)_N} $  $\Rightarrow$ 
 $ {\rm (ii)_N} $  $\Rightarrow$  ${\rm  (iii)_N} $ 
 $\Rightarrow$  ${\rm (ii)_{N+1}} $ and  ${\rm  (ii)_N} $  $\Rightarrow$ 
${\rm (iv)_N}:$

\noindent${\rm (i)_N}$~  $H^{naive}_n(A) =\{0\}\;$ for all $n \ge N$ and
 $H^{bar}_n(A) =\{0\}\;$ for all $n \ge N-1$;

\noindent${\rm (ii)_N}$~  $HH_n(A)= \{0\}\;$ for all $n \ge N;$ 

\noindent${\rm (iii)_N}$~  for all  $k \ge K$,
up to isomorphism of linear spaces,

$HC_{2k}(A) = HC_{N}(A)$ and
$HC_{2k+1}(A) = HC_{N-1}(A)$;

\noindent ${\rm (iv)_N}$~   up to isomorphism of linear spaces,
$HP_0(A) = HC_{N}(A)$
and 
$HP_1(A) = HC_{N-1}(A)$.\\
For Fr\'{e}chet algebras the isomorphisms in ${\rm (iii)_N}$ and
${\rm (iv)_N}$ are automatically topological. For a nuclear $\hat{\otimes}$-algebra ${\mathcal A}$ which is a  $DF$-space the isomorphisms in ${\rm (iii)_N}$
are automatically topological.
\end{proposition}

\begin{proof} A proof of the statement is given in \cite[Proposition 3.2]{Ly3}.
Here we add a part on the automatic continuity of the isomorphisms.
In view of the proofs of \cite[Propositions 2.1 and 3.2]{Ly3} 
 it is easy to see that isomorphisms of homology  in ${\rm (iii)_N}$ and
${\rm (iv)_N}$
are induced by continuous morphisms of complexes. 
Note that by Lemma \ref{topology_of_complexes},  for a  Fr\'{e}chet algebra  ${\mathcal A}$, the following complexes 
$(\mathcal{B}_c \bar{\Omega} {\mathcal A}_+, \tilde{b} +\tilde{B})$  and
$(\mathcal{B}_p \bar{\Omega} {\mathcal A}_+, \tilde{b} +\tilde{B})$ 
are complexes of Fr\'{e}chet spaces and continuous linear operators.
Thus, for \frechet\ algebras, by \cite[Lemma 0.5.9]{He0}, isomorphisms of homology groups are topological.

By Lemma \ref{topology_of_complexes}, for 
 a nuclear $\hat{\otimes}$-algebra ${\mathcal A}$ which is a  $DF$-space, the following complex 
$(\mathcal{B}_c \bar{\Omega} {\mathcal A}_+, \tilde{b} +\tilde{B})$ 
is a  complex of nuclear complete $DF$-spaces and continuous linear operators.
By Corollary \ref{OpenMap-strong-dual-nuclear}, for complete nuclear $DF$-spaces the isomorphisms for homology groups in ${\rm (iii)_N}$ are also topological.
\end{proof}

\begin{proposition}\label{n>N-HHn=0-cohomology}{\em \cite[Proposition 3.1]{Ly4}} Let $A$ be a  complete locally convex  algebra.
Then, for  any even integer $N$, say
$N=2K$,  and the following assertions, we have 
$ {\rm (i)_N} $  $\Rightarrow$ 
 $ {\rm (ii)_N} $  $\Rightarrow$  ${\rm  (iii)_N} $ 
 $\Rightarrow$  ${\rm (ii)_{N+1}} $ and  ${\rm  (ii)_N} $  $\Rightarrow$ 
${\rm (iv)_N}:$

\noindent ${\rm (i)_N}$~ 
$H_{naive}^n(A) =\{0\}\;$ for all $n \ge N$ and
 $H_{bar}^n(A) =\{0\}\;$ for all $n \ge N-1$;

\noindent ${\rm (ii)_N}$~   for all $n \ge N,$ 
$HH^n(A)= \{0\}$;

\noindent ${\rm (iii)_N}$~  for all  $k \ge K$,
up to isomorphism of linear spaces,
$HC^{2k}(A) = HC^{N}(A)$ and
$\;\;\;\;HC^{2k+1}(A) = HC^{N-1}(A)$;

\noindent ${\rm (iv)_N}$~   up to isomorphism of linear spaces,
$HP^0(A) = HC^{N}(A)$
and 
$HP^1(A) = HC^{N-1}(A)$.\\
For nuclear Fr\'{e}chet algebras the isomorphisms in ${\rm (iii)_N}$ and
${\rm (iv)_N}$ are topological isomorphisms. For a nuclear $\hat{\otimes}$-algebra ${\mathcal A}$ which is a  $DF$-space the isomorphisms in ${\rm (iii)_N}$ are topological isomorphisms.
\end{proposition}

\begin{proof} We need to add to the proof of \cite[Proposition 3.1]{Ly4}
the following part on automatic continuity.
In view of the proof of \cite[Proposition 3.1]{Ly4} it is easy to see that the isomorphisms of cohomology groups in ${\rm (iii)_N}$ and
${\rm (iv)_N}$ are induced by continuous morphisms of complexes. 

For nuclear \frechet\ algebras, by Lemma \ref{topology_of_complexes}, the 
complexes $((\mathcal{B}_c \bar{\Omega} {\mathcal A}_+)^*, \tilde{b}^* +\tilde{B}^*)$ 
and 
$((\mathcal{B}_p \bar{\Omega} {\mathcal A}_+)^*, \tilde{b}^* +\tilde{B}^*)$ are  complexes of strong duals of nuclear Fr\'{e}chet spaces.
By Lemma \ref{OpenMap-strong-dual},
the isomorphisms of cohomology groups in ${\rm (iii)_N}$ and
${\rm (iv)_N}$ are topological.

 For  a nuclear $\hat{\otimes}$-algebra ${\mathcal A}$ which is a  $DF$-space, by Lemma \ref{topology_of_complexes} and by \cite[Theorem 4.4.13]{Pi}, the chain  complex  
 $(\mathcal{B}_c \bar{\Omega} {\mathcal A}_+, \tilde{b} +\tilde{B})$ is the strong  dual of a complex of nuclear Fr\'{e}chet spaces. 
By \cite[Theorem 4.4.12]{Pi},  complete
nuclear $DF$-spaces and  nuclear Fr\'{e}chet spaces are reflexive.
Therefore $((\mathcal{B}_c \bar{\Omega} {\mathcal A}_+)^*, \tilde{b}^* +\tilde{B}^*)$ is a complex of nuclear Fr\'{e}chet spaces.
Thus, by \cite[Lemma 0.5.9]{He0}, the
isomorphisms of cohomology groups in ${\rm (iii)_N}$ are topological.
\end{proof}

The space of continuous traces on a topological algebra ${\mathcal A}$ is denoted by
${\mathcal A}^{tr}$, that is,
$${\mathcal A}^{tr} = \{f \in {\mathcal A}^*: f(ab) = f(ba)\;{\rm for\; all}\;a, b \in {\mathcal A} \}.$$
The closure in ${\mathcal A}$ of the linear span of elements of the form 
$\{ab-ba:\; a, b \in {\mathcal A} \}$ is denoted by $[{\mathcal A}, {\mathcal A}].$
Recall that $b_0: {\mathcal A} \hat{\otimes} {\mathcal A} \to {\mathcal A}$ is uniquely 
determined by $a
\otimes b \mapsto ab -ba$.

\begin{proposition}\label{A-simpl-trivial}
Let  ${\mathcal A}$ be in $\Phi$ and be
 a  $\hat{\otimes}$-algebra.

\noindent{\rm (i)} Suppose that the continuous cohomology groups 
${\mathcal H}^{naive}_{n}({\mathcal A}) = \{0 \}$ for  all  $n \ge 1$ and 
${\mathcal H}^{bar}_n({\mathcal A})= \{0\}\;$ for  all  $n \ge 0$. Then,
 up to  isomorphism of linear spaces,
\begin{equation}
\begin{array}{ccccccccccc} \label{nice-homology-lca}
{\mathcal H}{\mathcal H}_{n}({\mathcal A}) = \{0 \} \;\;{\rm for} \; 
{\rm  all} \;n \ge 1
\;{\rm and} \; {\mathcal H}{\mathcal H}_{0}({\mathcal A})  ={\mathcal A}/{\rm Im~} b_0;\\

{\mathcal H}{\mathcal C}_{2\ell}({\mathcal A}) = {\mathcal A}/{\rm Im~} b_0 \;{\rm and} \;
{\mathcal H}{\mathcal C}_{2\ell+1}({\mathcal A}) = \{0\} \;{\rm for} \; {\rm  all}\;
  \ell \ge 0; \\

{\mathcal H}{\mathcal P}_0({\mathcal A}) =  {\mathcal A}/{\rm Im~} b_0 \; {\rm and} \;
{\mathcal H}{\mathcal P}_1({\mathcal A}) = \{0\}.\\
\end{array}
\end{equation}


\noindent {\rm (ii)} Suppose that the continuous cohomology groups 
${\mathcal H}_{naive}^{n}({\mathcal A}) = \{0 \}$ for  all  $n \ge 1$ and 
${\mathcal H}_{bar}^n({\mathcal A})= \{0\}\;$ for  all  $n \ge 0$. Then,
up to  isomorphism of linear spaces,
\begin{equation}
\begin{array}{ccccccccccc} \label{nice-cohomology-lca}
{\mathcal H}{\mathcal H}^{n}({\mathcal A}) = \{0 \} \;\;{\rm for} \; {\rm  all}\;
 n \ge 1 \;{\rm and} \;
{\mathcal H}{\mathcal H}^{0}({\mathcal A})  ={\mathcal A}^{tr};\\

{\mathcal H}{\mathcal C}^{2\ell}({\mathcal A}) = {\mathcal A}^{tr} \;{\rm and} \;
\;\;\;\;{\mathcal H}{\mathcal C}^{2\ell+1}({\mathcal A}) = \{0\}  
\;\;{\rm for} \; {\rm  all}\;  \ell \ge 0; \\

{\mathcal H}{\mathcal P}^0({\mathcal A}) =  {\mathcal A}^{tr}
 \;{\rm and} \;
{\mathcal H}{\mathcal P}^1({\mathcal A}) = \{0\}.\\
\end{array}
\end{equation}
\end{proposition}

\begin{proof}  {\rm (i)}. One can see that 
${\mathcal H}^{bar}_n({\mathcal A})= \{0\}\; {\rm for \; all} \; n \ge 0\;$ implies that 
$${\mathcal H}{\mathcal H}_{n}({\mathcal A})=
{\mathcal H}^{naive}_{n}({\mathcal A})\; {\rm for \; all} \; n \ge 0,$$
see \cite[Section 3]{Ly3}. Note that by 
definition of  the `naive' Hochschild homology of ${\mathcal A}$, 
${\mathcal H}^{naive}_{0}({\mathcal A})  ={\mathcal A}/{\rm Im~} b_0$. Therefore,
${\mathcal H}{\mathcal H}_{n}({\mathcal A})= \{0\}$ for  all  $n \ge 1$  and
${\mathcal H}{\mathcal H}_{0}({\mathcal A})  ={\mathcal A}/{\rm Im~} b_0$.

From the exactness of the long Connes-Tsygan sequence of continuous
homology it follows that 
$${\mathcal H}{\mathcal C}_{0}({\mathcal A}) ={\mathcal H}_0^{naive}({\mathcal A}) ={\mathcal A}/{\rm Im~} b_0 
\;\;{\rm and}\;\; {\mathcal H}{\mathcal C}_{1}({\mathcal A}) =  \{ 0\}.$$
The rest of Statement (i) follows from  Proposition \ref{n>N-HHn=0-homology}. 

 {\rm (ii)} It is known that
${\mathcal H}_{bar}^n({\mathcal A})= \{0\}\; {\rm for \; all} \; n \ge 0,\;$
implies 
${\mathcal H}{\mathcal H}^{n}({\mathcal A})= 
{\mathcal H}_{naive}^{n}({\mathcal A})$ for all $ n \ge 0.$
By  definition of the `naive' Hochschild cohomology of ${\mathcal A}$, 
${\mathcal H}^0_{naive}({\mathcal A}) ={\mathcal A}^{tr}.$ 
Thus ${\mathcal H}{\mathcal H}^{n}({\mathcal A}) = \{0 \}$ for  all  $n \ge 1$ and
${\mathcal H}{\mathcal H}^{0}({\mathcal A})  ={\mathcal A}^{tr}$.

From the exactness of the long Connes-Tsygan sequence of continuous
 cohomology it follows that 
${\mathcal H}{\mathcal C}^{0}({\mathcal A}) ={\mathcal H}^0_{naive}({\mathcal A}) ={\mathcal A}^{tr}
\;\;{\rm and}\;\; {\mathcal H}{\mathcal C}^{1}({\mathcal A}) =  \{ 0\}.$
The rest of Statement (ii) follows from  Proposition \ref{n>N-HHn=0-cohomology}.
\end{proof}

\section{Cyclic-type cohomology of biflat 
$\hat{\otimes}$-algebras}

Recall that a  $\hat{\otimes}$-algebra ${\mathcal A}$ is said to be 
{\it biflat} if it is flat in
the category of  ${\mathcal A}$-$\hat{\otimes}$-bimodules  \cite[Def. 7.2.5]{He0}. A  $\hat{\otimes}$-algebra ${\mathcal A}$ is said to be 
{\it biprojective} if it is projective in
the category of  ${\mathcal A}$-$\hat{\otimes}$-bimodules  \cite[Def. 4.5.1]{He0}. By \cite[Proposition 4.5.6]{He0}, a 
$\hat{\otimes}$-algebra ${\mathcal A}$ is biprojective if and only if 
there exists an ${\mathcal A}$-$\hat{\otimes}$-bimodule morphism
$\rho_{{\mathcal A}}: {\mathcal A} \to {\mathcal A} \hat{\otimes} {\mathcal A}$ such that
$\pi_{{\mathcal A}} \circ \rho_{{\mathcal A}} = {\rm id}_{{\mathcal A}} $, where $\pi_{{\mathcal A}}$ is 
the canonical morphism $\pi_{{\mathcal A}}:{\mathcal A} \hat{\otimes}{\mathcal A} \to {\mathcal A} ,\;
 a_1\otimes a_2 \mapsto  a_1  a_2 .$
It can be proved that any biprojective $\hat{\otimes}$-algebra 
is biflat and ${\mathcal A} = \overline{{\mathcal A}^2} = {\rm Im~} \pi_{{\mathcal A}}$ \cite[Proposition 4.5.4]{He0}.
Here $\overline{{\mathcal A}^2}$ is the closure 
 of the linear span of the set
$\{a_1\cdot a_2:  a_1, a_2 \in {\mathcal A}\}$ in ${\mathcal A}$.
A  $\hat{\otimes}$-algebra ${\mathcal A}$ is said to be 
{\it contractible} if ${\mathcal A}_+$ is is projective in
the category of  ${\mathcal A}$-$\hat{\otimes}$-bimodules.
 A  $\hat{\otimes}$-algebra ${\mathcal A}$ is
contractible if and only if ${\mathcal A}$ is biprojective and has an identity
 \cite[Def. 4.5.8]{He0}.
For biflat Banach algebras  ${\mathcal A}$, Helemskii 
proved ${\mathcal A} = \overline{{\mathcal A}^2}= {\rm Im~} \pi_{{\mathcal A}}$ \cite[Proposition 7.2.6]{He0} and gave the description of the cyclic homology ${\mathcal H}{\mathcal C}_*$ and cohomology ${\mathcal H}{\mathcal C}^*$ groups of $A$ in \cite{He1}. Later the author generalized Helemskii's result to inverse limits of biflat Banach algebras \cite[Theorem 6.2]{Ly3} and to 
locally convex strict inductive limits of amenable Banach algebras \cite[Corollary 4.9]{Ly4}.

\begin{proposition}\label{biflat-A-trivial-lca}
Let  ${\mathcal A}$ be in $\Phi$ and be a biflat  
$\hat{\otimes}$-algebra such that ${\mathcal A} = {\rm Im~} \pi_{{\mathcal A}}$; 
in particular, let ${\mathcal A} \in \Phi$ be a biprojective 
$\hat{\otimes}$-algebra. Then

\noindent{\rm (i)}
${\mathcal H}^{naive}_{n}({\mathcal A}) = \{0 \}$ for  all  $n \ge 1$,
$H^{naive}_0({\mathcal A}) = {\mathcal A}/{\rm Im~} b_0$ and 
${\mathcal H}^{bar}_n({\mathcal A})= \{0\}\;$ for  all  $n \ge 0$;

\noindent{\rm (ii)} for the homology groups 
${\mathcal H}{\mathcal H}_*$,
${\mathcal H}{\mathcal C}_*$ and  ${\mathcal H}{\mathcal P}_*$
of ${\mathcal A}$
we have the isomorphisms of linear spaces {\rm (\ref{nice-homology-lca})}.

\noindent If, furthermore,  ${\mathcal A}$ is a \frechet\ space or 
 ${\mathcal A}$ is a nuclear $DF$-space, then 
$H^{naive}_0({\mathcal A}) = {\mathcal A}/[{\mathcal A}, {\mathcal A}]$ 
is  Hausdorff, and, for a biflat ${\mathcal A}$, 
${\mathcal A} = \overline{{\mathcal A}^2}$ implies that 
${\mathcal A} = {\rm Im~} \pi_{{\mathcal A}}$. 
\end{proposition}

\begin{proof} 
By \cite[Theorem 3.4.25]{He0}, up to topological isomorphism, the homology groups 
$$H^{naive}_n({\mathcal A}) = {\mathcal H}_n( {\mathcal A}, {\mathcal A})= 
{\rm Tor}^{{\mathcal A}^{e}}_n({\mathcal A}, {\mathcal A}_+) $$ 
for all $n \ge 0$. Since ${\mathcal A}$ is biflat, by \cite[Proposition 7.1.2]{He0},  $H^{naive}_n({\mathcal A}) =\{0\}\;$ for all $n \ge 1$. 

By \cite[Theorem 3.4.26]{He0}, up to topological isomorphism, the homology groups
$$ H^{bar}_n({\mathcal A}) = {\mathcal H}_{n+1}({\mathcal A},{\mathbf C}) = 
{\rm Tor}^{{\mathcal A}}_{n+1}({\mathbf C},{\mathbf C})$$
for all $n \ge 0$, where  ${\mathbf C}$ is the trivial ${\mathcal A}$-bimodule. Note that, for the trivial ${\mathcal A}$-bimodule ${\mathbf C}$, there is a flat resolution
$$0  \leftarrow {\mathbf C} \leftarrow {\mathcal A}_+ \leftarrow  
{\mathcal A} \leftarrow 0$$
in the category of  left or right  ${\mathcal A}$-$\hat{\otimes}$-modules.
By \cite[Theorem 3.4.28]{He0},
$ H^{bar}_n({\mathcal A}) ={\rm Tor}^{{\mathcal A}}_{n+1}({\mathbf C},{\mathbf C})= \{0\}\; $ for all $n \ge 1$. 
By assumption, 
${\mathcal A} = {\rm Im~} \pi_{{\mathcal A}}$, hence
 $H^{bar}_0({\mathcal A}) = {\mathcal A}/{\rm Im~} \pi_{{\mathcal A}} = \{0\}$.
Thus the conditions of Proposition \ref{A-simpl-trivial} (i) are satisfied.

In the categories of \frechet\ spaces and complete nuclear $DF$-spaces,
the open mapping theorem holds -- see Corollary \ref{open-map-strong-dual-FM} for $DF$-spaces. Thus, by \cite[Propositions 3.3.5 and 7.1.2]{He0}, up to topological isomorphism, 
$H^{naive}_0({\mathcal A}) = {\rm Tor}^{{\mathcal A}^{e}}_0({\mathcal A}, {\mathcal A}_+)$ is  Hausdorff.
Since ${\mathcal A}$ is biflat, by \cite[Proposition 7.1.2]{He0}, 
${\rm Tor}^{{\mathcal A}}_{0}({\mathbf C}, {\mathcal A})$ is also Hausdorff.
By \cite[Proposition 3.4.27]{He0},
$\overline{{\mathcal A}^2} = {\rm Im~} \pi_{{\mathcal A}}$.
\end{proof}

A  $\hat{\otimes}$-algebra ${\mathcal A}$ is {\it amenable} if
 ${\mathcal A}_+$ is a flat ${\mathcal A}$-$\hat{\otimes}$-bimodule.
For a Fr\'{e}chet algebra ${\mathcal A}$  amenability is equivalent to
the following: for all Fr\'{e}chet   ${\mathcal A}$-bimodules $X$, 
${\mathcal H}_{0}({\mathcal A}, X)$ is Hausdorff and 
 ${\mathcal H}_n({\mathcal A}, X)= \{0\}\;$ for all $n \ge 1$.
Recall that an amenable Banach algebra ${\mathcal A}$ is biflat and has a
bounded approximate identity \cite[Theorem VII.2.20]{He0}.

\begin{lemma}\label{amenable} Let  ${\mathcal A}$ be  an amenable
$\hat{\otimes}$-algebra which is 
 a \frechet\ space or a nuclear $DF$-space.  Then
${\mathcal H}^{naive}_{n}({\mathcal A}) = \{0 \}$ for  all  $n \ge 1$,
$H^{naive}_0({\mathcal A}) = {\mathcal A}/[{\mathcal A}, {\mathcal A}]$ and 
${\mathcal H}^{bar}_n({\mathcal A})= \{0\}\;$ for  all  $n \ge 0$.
\end{lemma}
\begin{proof} In the categories of \frechet\ spaces and complete nuclear $DF$-spaces, the open mapping theorem holds. Therefore, by \cite[Theorem III.4.25 and Proposition 7.1.2]{He0}, up to topological isomorphism,
 for the trivial ${\mathcal A}$-bimodule ${\mathbf C}$,
$$ H^{bar}_n({\mathcal A}) = {\mathcal H}_{n+1}({\mathcal A},{\mathbf C}) = 
{\rm Tor}_{n+1}^{{\mathcal A}^{e}}({\mathbf C},{\mathcal A}_+) =  \{0\}\;$$
for all $n \ge 0$;
$$H^{naive}_n({\mathcal A}) = {\mathcal H}_n( {\mathcal A}, {\mathcal A})= 
{\rm Tor}^{{\mathcal A}^{e}}_n({\mathcal A}, {\mathcal A}_+) = \{0\}\;$$
for all $n \ge 1$ and $H^{naive}_0({\mathcal A}) ={\rm Tor}^{{\mathcal A}^{e}}_0({\mathcal A},{\mathcal A}_+)$ is Hausdorff, that is, 
$H^{naive}_0({\mathcal A}) = {\mathcal A}/[{\mathcal A}, {\mathcal A}]$.
\end{proof}

\begin{theorem}\label{A-simpl-trivial-Fr-DF}
Let  ${\mathcal A}$ be  a  $\hat{\otimes}$-algebra which is 
 a \frechet\ space or a nuclear $DF$-space. 
 Suppose that the continuous homology groups 
${\mathcal H}^{naive}_{n}({\mathcal A}) = \{0 \}$ for  all  $n \ge 1$,
${\mathcal H}^{naive}_{0}({\mathcal A})$ is Hausdorff and 
${\mathcal H}^{bar}_n({\mathcal A})= \{0\}\;$ for  all  $n \ge 0$. 
In particular, asssume that ${\mathcal A}$ is a biflat algebra such 
that ${\mathcal A} = \overline{{\mathcal A}^2}$ or ${\mathcal A}$ is amenable.
Then

\noindent{\rm (i)} up to topological isomorphism, 
\begin{equation}
\begin{array}{ccccccccccc} \label{nice-homology}
{\mathcal H}{\mathcal H}_{n}({\mathcal A}) = \{0 \} \;\;{\rm for} \; {\rm  all}\;n \ge 1
\;{\rm and} \; {\mathcal H}{\mathcal H}_{0}({\mathcal A})  ={\mathcal A}/[{\mathcal A}, {\mathcal A}];\\
{\mathcal H}{\mathcal C}_{2\ell}({\mathcal A}) = {\mathcal A}/[{\mathcal A}, {\mathcal A}] \;{\rm and} \;
{\mathcal H}{\mathcal C}_{2\ell+1}({\mathcal A}) = \{0\} \;{\rm for} \; {\rm  all}\; \ell \ge 0; \\
\end{array}
\end{equation}
\noindent {\rm (ii)} up to topological isomorphism for \frechet\ algebras and up to isomorphism of linear spaces for nuclear $DF$-algebras,
\begin{equation}
\begin{array}{ccccccccccc} \label{nice-homology-HP}
{\mathcal H}{\mathcal P}_0({\mathcal A}) =  {\mathcal A}/[{\mathcal A}, {\mathcal A}] \;{\rm and} \;
{\mathcal H}{\mathcal P}_1({\mathcal A}) = \{0\};\\
\end{array}
\end{equation}
\noindent {\rm (iii)} 
\begin{equation}
\begin{array}{ccccccccccc} \label{nice-cohomology}
{\mathcal H}_{naive}^{n}({\mathcal A}) = \{0 \} \;\;{\rm for} \; {\rm  all}\; n \ge 1;\\
{\mathcal H}_{bar}^n({\mathcal A})= \{0\};\;{\rm for} \;{\rm  all}\; n \ge 0;\\
\end{array}
\end{equation}
\noindent {\rm (iv)} up to topological isomorphism for nuclear \frechet\ algebras and nuclear $DF$-algebras and up to isomorphism of linear spaces for  \frechet\ algebras, 
\begin{equation}
\begin{array}{ccccccccccc} \label{nice-cohomology-HC}
{\mathcal H}{\mathcal H}^{n}({\mathcal A}) = \{0 \} \;\;{\rm for} \; {\rm  all}\; n \ge 1 \;{\rm and} \;
{\mathcal H}{\mathcal H}^{0}({\mathcal A})  ={\mathcal A}^{tr};\\
{\mathcal H}{\mathcal C}^{2\ell}({\mathcal A}) = {\mathcal A}^{tr} \;{\rm and} \;
\;\;\;\;{\mathcal H}{\mathcal C}^{2\ell+1}({\mathcal A}) = \{0\}  
\;\;{\rm for} \; {\rm  all}\;  \ell \ge 0; \\
\end{array}
\end{equation}
\noindent {\rm (v)} up to topological isomorphism for nuclear \frechet\ algebras and up to isomorphism of linear spaces for  \frechet\ algebras and for nuclear $DF$-algebras,
\begin{equation}
\begin{array}{ccccccccccc} \label{nice-cohomology-HP}
{\mathcal H}{\mathcal P}^0({\mathcal A}) =  {\mathcal A}^{tr}
 \;{\rm and} \;
{\mathcal H}{\mathcal P}^1({\mathcal A}) = \{0\}.\\
\end{array}
\end{equation}
\end{theorem}

\begin{proof} In view of  Proposition \ref{biflat-A-trivial-lca}
 and Lemma \ref{amenable}, a biflat algebra ${\mathcal A}$ such 
that ${\mathcal A} = \overline{{\mathcal A}^2}$ and an amenable ${\mathcal A}$ 
satisfy the conditions of the theorem.

 By Proposition \ref{n>N-homology-dual-cohomology}, firstly,
${\mathcal H}^{bar}_n({\mathcal A})= \{0\}\; {\rm for \; all} \; n \ge 0\;$
if and only if 
${\mathcal H}_{bar}^n({\mathcal A})= \{0\}\; {\rm for \; all} \; n \ge 0;$ and,
secondly,
 ${\mathcal H}^n_{naive}({\mathcal A}) =\{0\}$ for all $n \ge 1$  if and only if 
${\mathcal H}_n^{naive}({\mathcal A}) =\{0\}$ for all $n \ge 1$ and
${\mathcal H}_0^{naive}({\mathcal A})$
is Hausdorff.

By Proposition \ref{A-simpl-trivial}, we have  isomorphisms of linear
spaces in {\rm (i)} -- {\rm (v)}.
In Propositions \ref{n>N-HHn=0-homology} and \ref{n>N-HHn=0-homology}
we show also when the above isomorphisms are automatically topological.
\end{proof}

\begin{remark} {\rm Recall that, for a biflat Banach algebra ${\mathcal A}$, 
 ${\rm db}_w {\mathcal A} \le 2$ \cite[Theorem 6]{Se1}. 
By \cite[Theorem 5.2]{Kh1}, for a Banach algebra ${\mathcal A}$ of a finite 
weak bidimension ${\rm db}_w {\mathcal A}$, we have isomorphisms
between the entire cyclic cohomology and the periodic cyclic cohomology of ${\mathcal A}$, $HE^0({\mathcal A}) = HP^0({\mathcal A}) = {\mathcal A}^{tr}$
and 
$HE^1({\mathcal A}) = HP^1({\mathcal A}) = \{0\}.$
The entire cyclic cohomology $HE^k({\mathcal A})$ of ${\mathcal A}$ for
$k=0,1$ are defined in \cite[IV.7]{Co2}.
In \cite[Theorem 6.1]{Pu} M. Puschnigg extended 
M. Khalkhali's result on the isomorphism  
$HE^k ({\mathcal A})= HP^k({\mathcal A})\;$ for $k=0,1$
from Banach algebras to some \frechet\ algebras.
}
\end{remark}


The following statement shows that the above theorems give the
explicit description of cyclic type homology and cohomology of 
the projective tensor product of two biprojective  
$\hat{\otimes}$-algebras.

\begin{proposition}\label{BtimesC-biprojective}
Let $\B$ and  ${\mathcal C}$ be biprojective 
 $\hat{\otimes}$-algebras. Then the projective tensor product 
${\mathcal A} = \B \hat{\otimes} {\mathcal C}$ is 
a biprojective  $\hat{\otimes}$-algebra.
\end{proposition}
\begin{proof} Since $\B$ is biprojective, there is a morphism
 of $\B$-$\hat{\otimes}$-bimodules 
$\rho_{\B}: \B \to \B \hat{\otimes} \B$ such that
$\pi_{\B} \circ \rho_{\B} = {\rm id}_{\B} $. A similar statement is
valid for ${\mathcal C}$. Let $i$ be the topological isomorphism
$${i}: 
(\B \hat{\otimes} \B) \hat{\otimes} ( {\mathcal C} \hat{\otimes} {\mathcal C}) \to
(\B \hat{\otimes} {\mathcal C}) \hat{\otimes} ( \B \hat{\otimes} {\mathcal C})$$
given by $(b_1\otimes b_2)\otimes (c_1\otimes c_2) \mapsto
(b_1\otimes c_1)\otimes (b_2\otimes c_2)$.
Note that $\pi_{\B \hat{\otimes} {\mathcal C}} = (\pi_{\B} \hat{\otimes}
\pi_{{\mathcal C}}) \circ {i}^{-1}$.
 It is routine to check that 
$$\rho_{\B \hat{\otimes} {\mathcal C}}:\B \hat{\otimes} {\mathcal C} \to 
(\B \hat{\otimes} {\mathcal C}) \hat{\otimes} ( \B \hat{\otimes} {\mathcal C})$$
defined by  $\rho_{\B \hat{\otimes} {\mathcal C}} = {i} \circ
(\rho_{\B} \otimes \rho_{{\mathcal C}})$
is a morphism of $\B \hat{\otimes} {\mathcal C}$-$\hat{\otimes}$-bimodules and
$\pi_{\B \hat{\otimes} {\mathcal C}} \circ \rho_{\B \hat{\otimes} {\mathcal C}} = 
{\rm id}_{\B \hat{\otimes} {\mathcal C}} $.
\end{proof}

\begin{remark} {\rm For amenable Banach algebras $\B$ and  ${\mathcal C}$,
B. E. Johnson showed that the Banach algebra 
 ${\mathcal A} = \B \hat{\otimes} {\mathcal C}$ is amenable \cite{BEJ1}.
By \cite[Proposition 5.4]{LW}, for a biflat Banach algebra ${\mathcal A}$,
each closed two-sided ideal $I$ with bounded approximate identity is
amenable and the quotient algebra ${\mathcal A}/I$ is biflat. 
Thus the explicit description of cyclic type homology and cohomology
of such $I$ and ${\mathcal A}/I$ is also given in Theorem 
\ref{A-simpl-trivial-Fr-DF}.
One can find a number of examples of biflat and simplicially trivial
Banach and $C^*$- algebras in \cite[Example 4.6, 4.9]{Ly3}.
}
\end{remark}

\section{Applications to the cyclic-type cohomology of biprojective 
$\hat{\otimes}$-algebras}
In this section we present  examples of  nuclear  biprojective 
$\hat{\otimes}$-algebras which are \frechet\ spaces or 
$DF$-spaces and the continuous cyclic-type homology and cohomology
of these algebras.

\begin{example}~ 
{\rm Let  $G$ be a compact Lie group and let 
$\mathcal{E}(G)$ be the nuclear \frechet\ algebra  of 
smooth functions on $G$  with the convolution product.
It was shown by Yu.V. Selivanov that ${\mathcal A} =\mathcal{E}(G)$ is 
biprojective \cite{Se3}.

Let $\mathcal{E}^*(G)$ be the strong dual to 
$\mathcal{E}(G)$, so that  $\mathcal{E}^*(G)$ is a complete nuclear $DF$-space. 
This is a $\hat{\otimes}$-algebra with respect to convolution 
multiplication: 
for $f,g \in \mathcal{E}^*(G)$ and $x \in  \mathcal{E}(G)$, 
$ <f*g, x>= <f, y>$, where $y \in \mathcal{E}(G)$ is defined by 
$ y(s)= <g, x_s>, \; s \in G$ and $x_s(t)= x (s^{-1}t),\; t \in G$.
J.L. Taylor proved that the algebra of distributions
$\mathcal{E}^*(G)$ on a compact Lie group $G$ is contractible 
 \cite{Tay1}.
}
\end{example}
\begin{example}~ {\rm  Let $(E, F)$ be a pair of complete Hausdorff
locally convex spaces endowed with a jointly continuous bilinear 
form $ \langle \cdot, \cdot \rangle : E \times F \to {\mathbf C}$
that is not identically zero. The space ${\mathcal A}=E \hat{\otimes} F $ is
a $\hat{\otimes}$-algebra with respect to the multiplication defined
by 
$$ (x_1 \otimes x_2) (y_1 \otimes y_2) = \langle x_2, y_1 \rangle  
x_1 \otimes y_2, \; x_i \in E,\; y_i \in F.$$
Yu.V. Selivanov proved that this algebra is biprojective and
 usually non unital \cite{Se2,Se3}. More exactly, if 
${\mathcal A}=E \hat{\otimes} F $ has a left or right identity, then
$E$ or $F$ respectively is finite-dimensional. If the form 
$\langle \cdot, \cdot \rangle$
is nondegenerate, then ${\mathcal A}=E \hat{\otimes} F $ is called the {\it
tensor algebra generated by the duality} 
$(E, F, \langle \cdot, \cdot \rangle)$.

In particular, if $E$ is a Banach space with the approximation
property, then the algebra ${\mathcal A}=E \hat{\otimes} E^* $ is isomorphic to
the algebra ${\mathcal N}(E)$ of nuclear operators on $E$
\cite[II.2.5]{He0}.
}
\end{example}

\subsection{ K\"{o}the  sequence algebras}~ 
The following results on K\"{o}the algebras can be found in A. Yu.
Pirkovskii's papers \cite{Pir1,Pir2}.

 A set $P$ of nonnegative real-valued sequences 
$p = (p_i)_{i \in \mathbf{N}}$ is called  a {\it K\"{o}the set} if 
the following axioms are satisfied:

($P1$) for every $i \in \mathbf{N}$ there is $p \in P$ such that $p_i
>0$;

($P2$) for every $p, q \in P$ there is $r \in P$ such that
${\rm max} \{p_i, q_i\} \le r_i$ for all $i \in \mathbf{N}$.

\noindent Suppose, in addition, the following condition is satisfied:
 
($P3$) for every $p \in P$ there exist  $q \in P$ and a constant
$C>0$ such that
$p_i \le C q_i^2$ for all $i \in \mathbf{N}$.

For any K\"{o}the set $P$ which satisfies ($P3$),
the K\"{o}the space
$$ \lambda(P)= \{ x = (x_n) \in {\mathbf C}^\mathbf{N} : \|x \|_p = \sum_n |x_n| p_n
< \infty \;{\rm for} \;{\rm all} \; p \in P \}$$
is a complete locally convex space with the topology determined by
the family of seminorms  $\{\|x \|_p :p \in P \}$ and
 a $\hat{\otimes}$-algebra with pointwise multiplication.
The $\hat{\otimes}$-algebras $ \lambda(P)$ are called 
 {\it K\"{o}the algebras}.

\noindent By \cite{Pi} and \cite{GrA}, for a  K\"{o}the set, 
 $ \lambda(P)$ is  {\it nuclear} if and only if  

($P4$) for every $p \in P$ 
there  exist $q \in P$ and $\xi \in \ell^1$ 
 such that $p_i \le \xi_i q_i$ for all $i \in \mathbf{N}$.

\noindent By \cite[Theorem 3.5]{Pir1}, $ \lambda(P)$ is {\it biprojective} if and
only if 

($P5$) for every $p \in P$ there exist $q \in P$ and a constant
$M>0$ such that $p_i^2 \le M q_i$ for all $i \in \mathbf{N}$.

The algebra $ \lambda(P)$ is {\it unital} if and only if 
$\sum_n p_n < \infty$  for every $p \in P$.

\begin{example}~{\rm  Fix a real number $1 \le R \le \infty$ and a
nondecreasing sequence $\alpha = (\alpha_i)$ of positive numbers with
$\lim_{i \to  \infty} \alpha_i = \infty$. The power series space
$$ \Lambda_R(\alpha)= 
\{ x = (x_n) \in {\mathbf C}^\mathbf{N} : \|x \|_r = \sum_n |x_n| r^{\alpha_n}
< \infty \;{\rm for} \;{\rm all} \; 0 <r <R \}$$
is a \frechet\  K\"{o}the algebra with pointwise multiplication.
The topology  of $ \Lambda_R(\alpha)$ is  determined by
a countable  family of seminorms  $\{\|x \|_{r_k} :k \in \mathbf{N}\}$ 
where $\{r_k\}$ is an arbitrary increasing sequence converging to $R$.

By \cite[Corollary 3.3]{Pir2}, $ \Lambda_R(\alpha)$ is {\it biprojective }
 if and only if $R=1$ or $R=\infty$.

By the Grothendieck-Pietsch criterion, $\Lambda_R(\alpha)$ is {\it nuclear} if and only if  for $ \overline{\lim_n }\frac{\log n}{\alpha_n} = 0$ for
$R < \infty$ and $ \overline{\lim_n }\frac{\log n}{\alpha_n}  < \infty$
for $R = \infty$, see  \cite[Example 3.4]{Pir1}.
}
\end{example}

The algebra $ \Lambda_R((n))$ is topologically isomorphic to the 
algebra of functions holomorphic on the open disc of radius $R$,
 endowed with {\it Hadamard product}, that is, with 
``co-ordinatewise" product of the Taylor expansions of holomorphic 
functions.

\begin{example}~{\rm
The algebra $\mathcal{H}({\mathbf C}) \iso  \Lambda_{\infty}((n))$ of entire functions, endowed with the Hadamard product, is a
biprojective nuclear \frechet\ algebra \cite{Pir2}.
}
\end{example}
\begin{example}~{\rm
The algebra $\mathcal{H}({\mathbf D}_1) \iso \Lambda_1((n))$ of functions holomorphic on the open unit disc,
 endowed with the Hadamard product, is a
biprojective nuclear \frechet\ algebra. Moreover it is contractible,
since the function $z \mapsto (1- z)^{-1}$ is an identity for
 $\mathcal{H}({\mathbf D}_1)$ \cite{Pir2}.
}
\end{example}

For any  K\"{o}the space $ \lambda(P)$ the dual space $\lambda(P)^*$
can be canonically identified with 
$$ \{(y_n) \in {\mathbf C}^\mathbf{N} : \exists
 p \in P \;  {\rm and} \; C >0
\;{\rm such} \;{\rm that} \;
|y_n| \le C p_n \;{\rm for} \;{\rm all} \; n \in \mathbf{N} \}.$$
It is shown in \cite{Pir2} that, for a biprojective
K\"{o}the algebra $ \lambda(P)$, $ \lambda(P)^*$ is a sequence algebra
with pointwise multiplication.

The algebra $ \lambda(P)^*$ is {\it unital} if and only if there exists 
 $p \in P$ such that $\inf_i p_i >0$.

\begin{example}\label{s-algebra}~{\rm
 The nuclear \frechet\ algebra of rapidly decreasing
sequences 
 $$s = \{ x = (x_n) \in {\mathbf C}^\mathbf{N} : \|x \|_k = \sum_n |x_n| n^k
< \infty \;{\rm for} \;{\rm all} \; k \in \mathbf{N} \}$$
is a biprojective  K\"{o}the algebra \cite{Pir1}. The algebra $s$
is topologically isomorphic to $ \Lambda_{\infty}(\alpha)$
with $\alpha_n =\log n$ \cite{Pir2}. The nuclear K\"{o}the 
$\hat{\otimes}$-algebra $s^*$ 
of sequences of polynomial growth is contractible \cite{Tay1}.
}
\end{example} 

\begin{example} \cite[Section 4.2]{Pir2}~ {\rm 
Let $P$ be a K\"{o}the set such that $p_i \ge 1$ for all $p \in P$ and all
$n \in \mathbf{N}$. Then the formula $  \langle a,b \rangle  = \sum_i a_i b_i $ 
defines a jointly continuous, nondegenerate bilinear form on $\lambda (P)
\times \lambda (P)$. Thus  $M(P)=\lambda (P) \hat{\otimes}\lambda (P) $ can be
considered as  the tensor algebra generated by the duality 
$(\lambda (P), \lambda (P), \langle \cdot, \cdot \rangle)$, and so 
is biprojective.
There is a canonical isomorphism between $M(P)$ and the algebra 
$\lambda (P \times P)$ of
$\mathbf{N} \times \mathbf{N}$ complex matrices
$(a_{ij})_{(ij) \in \mathbf{N} \times \mathbf{N}}$ 
satisfying the condition
$ \| a \|_p = \sum_{i,j} |a_{ij}| p_i p_j < \infty$
for all $p \in P$ with the usual matrix multiplication.

In particular, for $P=\{ (n^k)_{n \in  \mathbf{N}}: k = 0,1, \dots \}$,
we obtain the biprojective nuclear \frechet\
algebra $\Re = s \hat{\otimes} s$ of ``smooth compact
operators" consisting of $\mathbf{N} \times \mathbf{N}$ complex matrices
$(a_{ij})$  with rapidly decreasing matrix entries. Here $s$ is from
Example \ref{s-algebra}. 
}
\end{example}
\begin{theorem}\label{cyclic-biproj-examples}
 Let ${\mathcal A}$ be a $\hat{\otimes}$-algebra belonging to
one of the following classes:

{\rm (i)} ${\mathcal A} =\mathcal{E}(G)$ or
 ${\mathcal A} =\mathcal{E}^*(G)$  for a compact Lie group $G$;

{\rm (ii)} ${\mathcal A}=E \hat{\otimes} F$, the 
tensor algebra generated by the duality  $(E, F, \langle \cdot, \cdot
\rangle)$ for nuclear \frechet\ spaces $E$ and $F$ 
(e.g., $\Re = s \hat{\otimes} s$) or for 
nuclear complete  $DF$-spaces $E$ and $F$;

{\rm (iii)}  \frechet\ K\"{o}the algebras
  ${\mathcal A} =\lambda(P)$ such that  the  K\"{o}the set $P$ 
satisfies {\rm ($P3$), ($P4$)} and {\rm ($P5$)}; 
in particular, $ \Lambda_1(\alpha)$ such that 
 $ \overline{\lim_n }\frac{\log n}{\alpha_n} = 0$ or
$ \Lambda_{\infty}(\alpha)$ such that 
 $ \overline{\lim_n }\frac{\log n}{\alpha_n}  < \infty$.
(e.g., $\mathcal{H}({\bf D}_1)$, $s$, $\mathcal{H}({\mathbf C})$).

{\rm (iv)}  K\"{o}the algebras 
${\mathcal A} = \lambda(P)^*$ which are the strong duals of $\lambda(P)$
from {\rm (iii)}.

{\rm (v)} the projective tensor product ${\mathcal A} = \B \hat{\otimes} {\mathcal C}$ of 
biprojective nuclear $\hat{\otimes}$-algebras $\B$ and  ${\mathcal C}$ which are 
\frechet\ spaces or $DF$-spaces; in particular, 
${\mathcal A}= \mathcal{E}(G)\hat{\otimes} \Re$.

Then, up to topological isomorphism,
$${\mathcal H}^{naive}_{n}({\mathcal A}) = \{0 \} \;\;{\rm for} \; 
{\rm  all}\; n \ge 1 \; {\rm and}\; {\mathcal H}^{naive}_0({\mathcal A}) = 
{\mathcal A}/[{\mathcal A}, {\mathcal A}];$$ 
$${\mathcal H}^{bar}_n({\mathcal A})= \{0\}\;\;{\rm for} \; 
{\rm  all}\;n \ge 0;$$ 
$${\mathcal H}{\mathcal H}_{n}({\mathcal A}) = \{0 \} \;\;{\rm for} \; 
{\rm  all}\;n \ge 1 \;{\rm and} \; {\mathcal H}{\mathcal H}_{0}({\mathcal A})  
={\mathcal A}/[{\mathcal A}, {\mathcal A}];$$ 
$${\mathcal H}{\mathcal C}_{2\ell}({\mathcal A}) = 
{\mathcal A}/[{\mathcal A}, {\mathcal A}] \;{\rm and} \;
{\mathcal H}{\mathcal C}_{2\ell+1}({\mathcal A}) = \{0\} \;{\rm for} \; {\rm  all}\;
  \ell \ge 0; $$ 
$${\mathcal H}_{naive}^{n}({\mathcal A}) = \{0 \} \;\;{\rm for} \; 
{\rm  all}\; n \ge 1; \;
{\mathcal H}_{bar}^n({\mathcal A})= \{0\}\;\;{\rm for} \; 
{\rm  all}\;n \ge 0;$$ 
$${\mathcal H}{\mathcal H}^{n}({\mathcal A}) = \{0 \} \;\;{\rm for} \; {\rm  all}\;
 n \ge 1 \;{\rm and} \;
{\mathcal H}{\mathcal H}^{0}({\mathcal A})  ={\mathcal A}^{tr};$$ 
$${\mathcal H}{\mathcal C}^{2\ell}({\mathcal A}) = {\mathcal A}^{tr} \;{\rm and} \;
\;\;\;\;{\mathcal H}{\mathcal C}^{2\ell+1}({\mathcal A}) = \{0\}  
\;\;{\rm for} \; {\rm  all}\;  \ell \ge 0;$$
and, up to topological isomorphism for \frechet\ algebras and up to
isomorphism of linear spaces for $DF$-algebras,
$${\mathcal H}{\mathcal P}_0({\mathcal A}) =  {\mathcal A}/[{\mathcal A}, {\mathcal A}] \;{\rm and} \;
{\mathcal H}{\mathcal P}_1({\mathcal A}) = \{0\};$$
$${\mathcal H}{\mathcal P}^0({\mathcal A}) =  {\mathcal A}^{tr}
 \;{\rm and} \;
{\mathcal H}{\mathcal P}^1({\mathcal A}) = \{0\}.$$
\end{theorem}
\begin{proof} We have mentioned above that the algebras in  
(i)-(iii) and (v) are biprojective and nuclear. 
By \cite[Corollary 3.10]{Pir2}, for any nuclear biprojective \frechet\
K\"{o}the algebra $ \lambda(P)$, the strong dual 
$ \lambda(P)^*$ is a nuclear, biprojective K\"{o}the 
$\hat{\otimes}$-algebra which is a  $DF$-space. For nuclear 
\frechet\ algebras and for nuclear $DF$-algebras,  the conditions of 
 Theorem \ref{A-simpl-trivial-Fr-DF} are satisfied. Therefore,
for the homology and cohomology groups ${\mathcal H}{\mathcal H}$ and
${\mathcal H}{\mathcal C}$ of ${\mathcal A}$
we have the topological isomorphisms
{\rm (\ref{nice-homology})} and {\rm (\ref{nice-cohomology-HC})}.
For the periodic cyclic homology and cohomology groups
 ${\mathcal H}{\mathcal P}$ of ${\mathcal A}$, for \frechet\ algebras, we have 
 topological isomorphisms and,  for nuclear $DF$-algebras,
isomorphisms of linear spaces {\rm (\ref{nice-homology-HP})} and {\rm (\ref{nice-cohomology-HP})}.
It is obvious that, for commutative algebras, ${\mathcal A}^{tr}=
{\mathcal A}^*$ and ${\mathcal A}/[{\mathcal A}, {\mathcal A}] ={\mathcal
A}$.
\end{proof}

The cyclic-type homology and cohomology of 
$\mathcal{E}(G)$  for a compact Lie group $G$ were calculated in \cite{Mey}.


\end{document}